\documentclass[twoside,twocolumn]{article}
\usepackage[twocolumn,textwidth=18cm,columnsep=8mm]{geometry}
\usepackage{graphicx,amsmath,amssymb}

\title{What is Aperiodic Order?}

\author{Michael Baake, David Damanik and Uwe Grimm}

\begin{document}

\maketitle

Crystals with their regular shape and pronounced faceting have
fascinated humans for ages. About a century ago, the understanding of
the internal structure increased considerably through the work of Max
von Laue (Nobel prize in physics 1914) and William Henry Bragg and
William Lawrence Bragg (father and son, joint Nobel prize in physics
1915). They developed X-ray crystallography and used it to show that a
lattice-periodic array of atoms lies at the heart of the matter.  This
interpretation became the accepted model for solids with pure Bragg
diffraction, which was later extended to allow for incommensurately
modulated structures.

\begin{figure}[ht]
\centerline{\includegraphics[width=0.95\columnwidth]{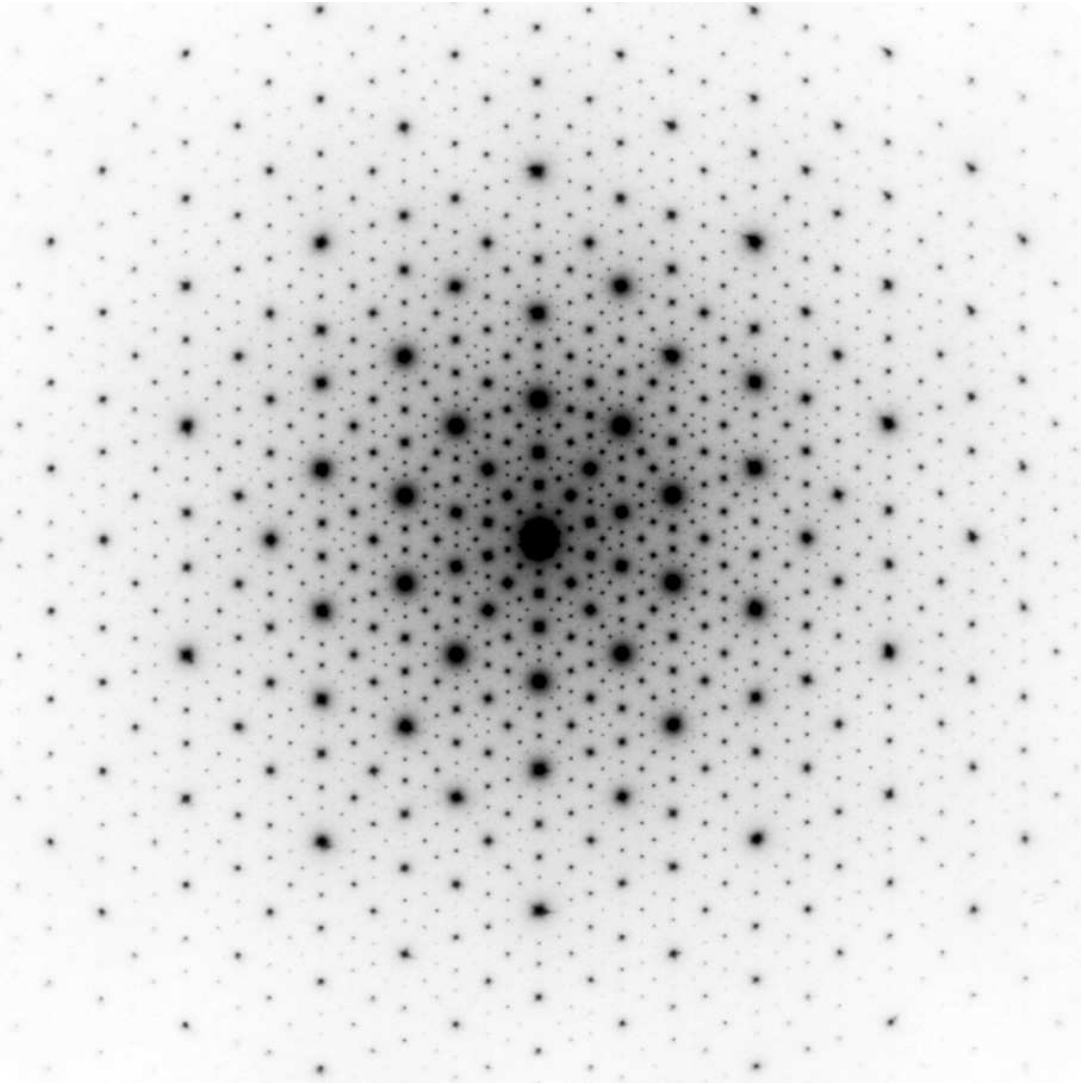}}
\caption{Electron diffraction image (intensity inverted) of an
  icosahedrally symmetric AlMnPd allow, taken along a fivefold
  axis. Due to the projection, the images is tenfold symmetric. Image
  courtesy of Conradin Beeli.\label{fig:ico}}
\end{figure}

In 1982, the materials scientist Dan Shechtman discovered a perfectly
diffractive solid with a non-crystallographic (icosahedral) symmetry
\cite{SBGC}; see Figure~\ref{fig:ico} for a qualitatively similar
experimental diffraction image.  This discovery, for which he received
the 2011 Nobel prize in chemistry, was initially met with disbelief
and heavy criticism, although such structures could have been expected
on the basis of Harald Bohr's work on almost periodic functions. In
fact, the situation is a classic case of a `missed opportunity'. Let
us try to illustrate this a little further and thus explain some
facets of what is now known as the theory of aperiodic order.

\begin{figure}[t]
\centerline{\includegraphics[width=0.95\columnwidth]{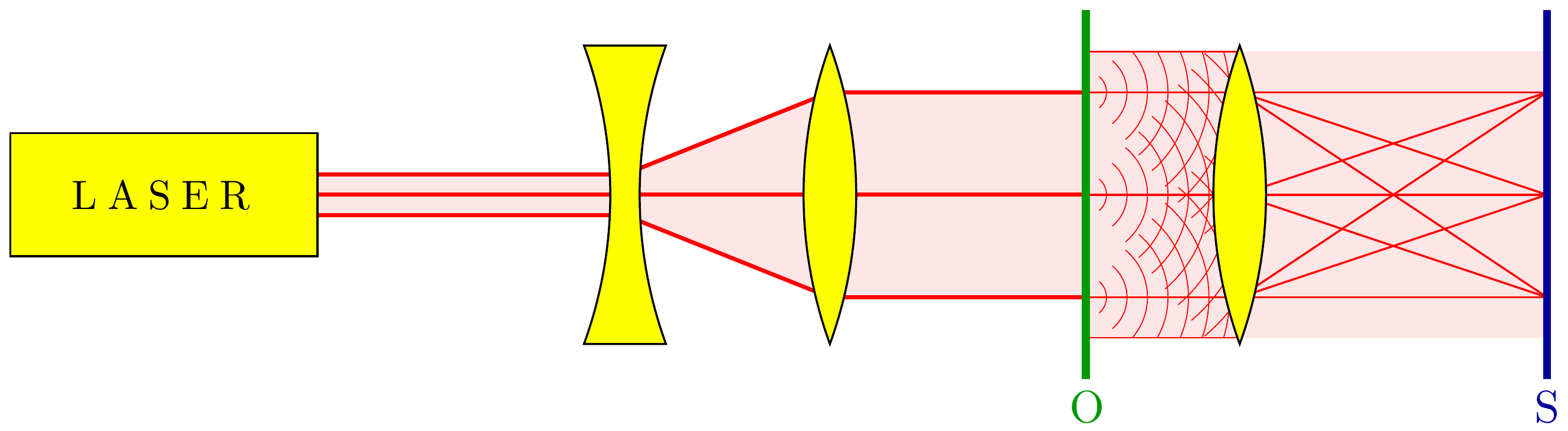}}
\caption{Schematic representation of a simple optical diffraction experiment.
The object (O) is illuminated by  a coherent light source,
and the diffracted intensity is collected on a screen (S), with sharp 
intensity peaks arising from  constructive interference of scattered waves.
\label{fig:diff}}
\end{figure}

Let us consider a uniformly discrete point set $\varLambda$ in
Euclidean space, where the points are viewed as idealizations of atomic
positions (all assumed to be of one type for simplicity).  Much of the
terminology for such point sets was developed by Jeffrey C.\ Lagarias
\cite{Lag}.  Placing unit point measures at each position in
$\varLambda$ leads to the associated Dirac comb
\[
    \delta^{}_{\!\varLambda} \, = \, \sum_{x\in\varLambda} \delta_{x}.
\]    
A diffraction experiment measures the correlation between atomic
locations.  Mathematically, this is expressed through the diffraction
measure, which is the Fourier transform
$\widehat{\,\gamma^{}_{\!\varLambda}\,}$ of the autocorrelation
\[
    \gamma^{}_{\!\varLambda} \, = \, \lim_{R\to\infty} \frac{
    \delta^{}_{\!\varLambda\cap B^{}_{\! R}} * 
    \delta^{}_{-\varLambda\cap B^{}_{\! R}}}{\mathrm{vol}(B^{}_{\! R})} ,
\]
provided that this limit exists in the vague topology; see Hof's
contribution to \cite{Nato} or \cite{BG} for details. Here,
$\widehat{\,\gamma^{}_{\!\varLambda}\,}$ describes the outcome of a
(kinematic) diffraction experiment, such as the one sketched in
Figure~\ref{fig:diff} with an optical bench, which should be available
in most physics laboratories for experimentation.

For a lattice periodic point set, the diffraction measure is a pure
point measure supported on the dual lattice. This implies that a
tenfold symmetric diffraction diagram, such as the one of
Figure~\ref{fig:ico}, cannot be produced by a lattice periodic
structure, as lattices in two or three dimensions can only have two-,
three-, four- or sixfold rotational symmetry by the crystallographic
restriction. This raises the question what types of point sets can
generate such new kinds of diffraction measures, which are pure point
and display non-crystallographic symmetries.

\begin{figure}[ht]
\centerline{\includegraphics[width=0.95\columnwidth]{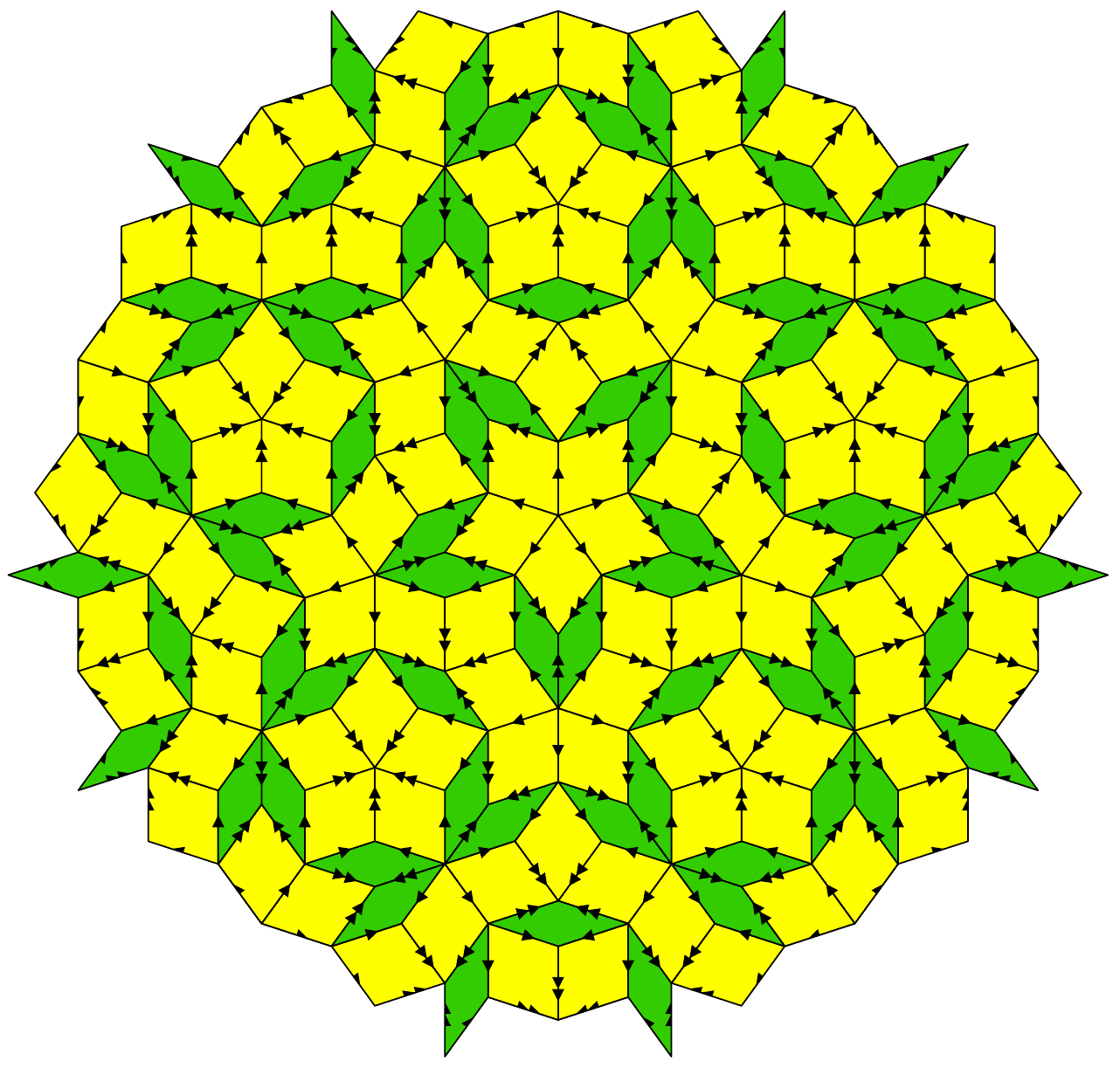}}
\caption{Fivefold symmetric patch of the rhombic Penrose tiling, 
which is equivalent to various other versions; see \cite[Sec.~6.2]{BG} 
for details. The arrow markings represent local rules that enforce 
aperiodicity.\label{fig:pen}}
\end{figure}

\begin{figure}[ht]
\centerline{\includegraphics[width=0.95\columnwidth]{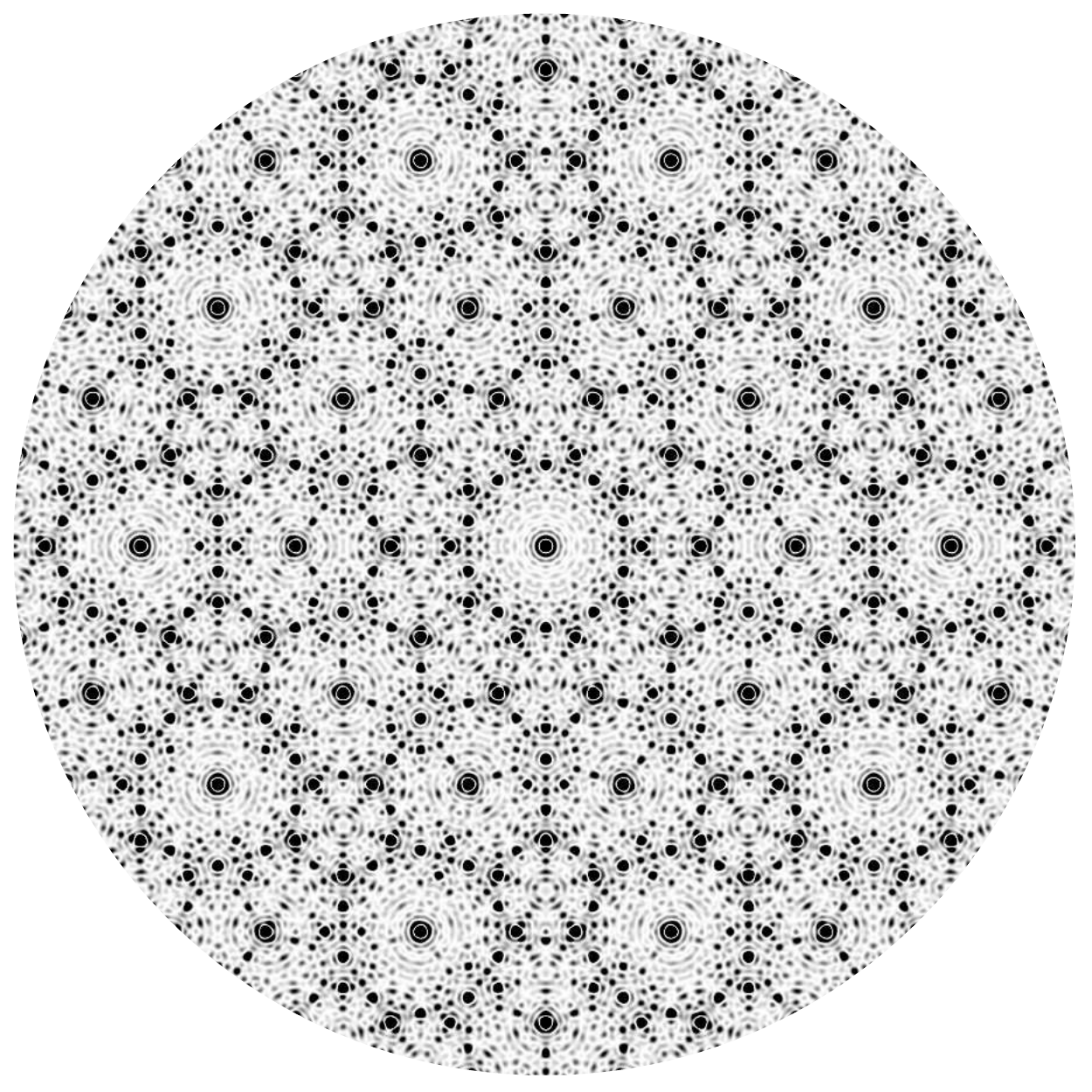}}
\caption{Circular detail of the (intensity inverted) diffraction image
  of the finite patch of Figure~\ref{fig:pen}, with point measures
  placed on all vertex points of the rhombic
  tiling.\label{fig:pendiff}}
\end{figure}

Let us begin with an example. Already in 1974, Roger Penrose
constructed an aperiodic tiling of the entire plane, equivalent to the
one shown in Figure~\ref{fig:pen}.  Taking this finite patch and
considering the set of vertices as the point set $\varLambda$, the
diffraction measure is the one shown in Figure~\ref{fig:pendiff},
which resembles what Alan L.\ Mackay observed when he performed an
optical diffraction experiment with an assembly of small disks
centered at the vertex points of a rhombic Penrose tiling, shortly
before the discovery of quasicrystals \cite{Mac}. For an infinite
tiling, the diffraction measure is pure point and tenfold symmetric.

The Penrose tiling is a particularly prominent example of a large
class of point sets, for which it has been shown that the diffraction
measure exists and is pure point. Such sets, which were first
introduced by Yves Meyer \cite{Mey} and are nowadays called model
sets, are constructed from a cut and project scheme (CPS)
\[
\renewcommand{\arraystretch}{1.2}\begin{array}{r@{}ccccc@{}l}
   & \mathbb{R}^{d} & \xleftarrow{\,\;\;\pi\;\;\,} & \mathbb{R}^{d} \times H & 
        \xrightarrow{\;\pi^{}_{\mathrm{int}\;}} & H & \\
   & \cup & & \cup & & \cup & \hspace*{-1ex} 
   \raisebox{1pt}{\text{\footnotesize dense}} \\
   & \pi(\mathcal{L}) & \xleftarrow{\; 1-1 \;} & \mathcal{L} & 
   \xrightarrow{\; \hphantom{1-1} \;} & \pi^{}_{\mathrm{int}}(\mathcal{L}) & \\
   & \| & & & & \| & \\
   & L & \multicolumn{3}{c}{\xrightarrow{\qquad\qquad\;\;\,\star\,
       \;\;\qquad\qquad}} 
       &  {L_{}}^{\star} & \\
\end{array}\renewcommand{\arraystretch}{1}
\]
where $\mathcal{L}$ is a lattice in $\mathbb{R}^d\times H$ whose
projection $\pi^{}_{\mathrm{int}}(\mathcal{L})$ to `internal space'
$H$ is dense. The projection $\pi$ to `physical space'
$\mathbb{R}^{d}$ is required to be injective on $\mathcal{L}$, so that
$\pi$ is a bijection between $\mathcal{L}$ and
$L=\pi(\mathcal{L})$. The CPS thus provides a well-defined mapping
$\star\! :\; L\to H$ with $x\mapsto
x^{\star}=\pi^{}_{\mathrm{int}}\bigl((\pi|_{\mathcal{L}})^{-1}(x)\bigr)$. The
internal space $H$ is often a Euclidean space, but the general theory
works for locally compact Abelian groups. A model set is then obtained
by choosing a `window' $W\subset H$ and defining the point set
\[
   \varLambda \, = \, \{ x\in L \mid x^{\star}\in W\} 
              \, \subset\,  \mathbb{R}^{d}.
\]
The Penrose point set arises in this framework with a two-dimensional
Euclidean internal space and the root lattice $A_{4}$, while the
generalisation to model sets with icosahedral symmetry was first
discussed by Peter Kramer \cite{KN}. It was Robert V.\ Moody
\cite{Moo,Nato} who recognized the connections to Meyer's abstract
concepts and championed their application in the theory of aperiodic
order and their further development into the shape and form used
today.

For a non-empty, compact window that is the closure of its interior
and whose boundary has Haar measure $0$, the resulting model set has a
diffraction measure that is supported on the projection
$\pi(\mathcal{L^{*}})$ of the dual lattice, and hence is a pure point
measure; see Schlottmann's article in \cite{BM} for a proof. The
diffraction intensities can be calculated explicitly in this setting.

While a pure point diffraction detects order, it is not true that all
ordered structures have pure point diffraction. A simple example can
be generated by the Thue--Morse substitution $0\mapsto 01$, $1\mapsto
10$, via taking a bi-infinite sequence that is invariant under
the square of this rule
\[
    \ldots 0110100110010110|0110100110010110 \ldots
\]
and taking only those $x\in\mathbb{Z}$ that correspond to the
positions of $1$s in the sequence. The diffraction measure of this
one-dimensional point set is not pure point, as it contains a
non-trivial singular continuous component.

There is a very useful connection between the diffraction measure and
the spectral measures of an associated dynamical system. Starting from
a point set $\varLambda\subset\mathbb{R}^{d}$ of finite local
complexity (which means that, up to translations, there are only
finitely many patches for any given size), we define its hull
\[
   \mathbb{X}(\varLambda) \, =  \, 
    \overline{\{x+\varLambda\mid x\in\mathbb{R}^{d}\}}
\]
where the closure is taken in the local topology (closeness forces
coincidence on a large ball around the origin up to a small
translation). $\mathbb{R}^{d}$ acts on $\mathbb{X}(\varLambda)$ by
translations. Assume that there is an invariant probability
measure $\mu$. This induces an action of $\mathbb{R}^{d}$ on
$L^{2}(\mathbb{X},\mu)$ via unitary operators. It is a fundamental
result that the diffraction measure is pure point if and only if all
spectral measures are pure point.

More generally, the spectral measures correspond to diffraction
measures associated with elements of certain topological factors of
the dynamical system $(\mathbb{X},\mathbb{R}^{d},\mu)$
\cite{BLvE}. For example, the diffraction measure of the Thue--Morse
point set consists of the pure point part $\delta^{}_{\mathbb{Z}}$ and
a non-trivial singular continuous component. The non-trivial point
part of the dynamical spectrum, which is $\mathbb{Z}[\frac{1}{2}]$, is
therefore not fully detected by the diffraction measure of the
Thue--Morse point set, and only shows up in the diffraction measure of
a global 2-to-1 factor, which is a model set, with $H=\mathbb{Z}_{2}$,
the $2$-adic numbers.

\begin{figure}[ht]
\centerline{\includegraphics[width=0.8\columnwidth]{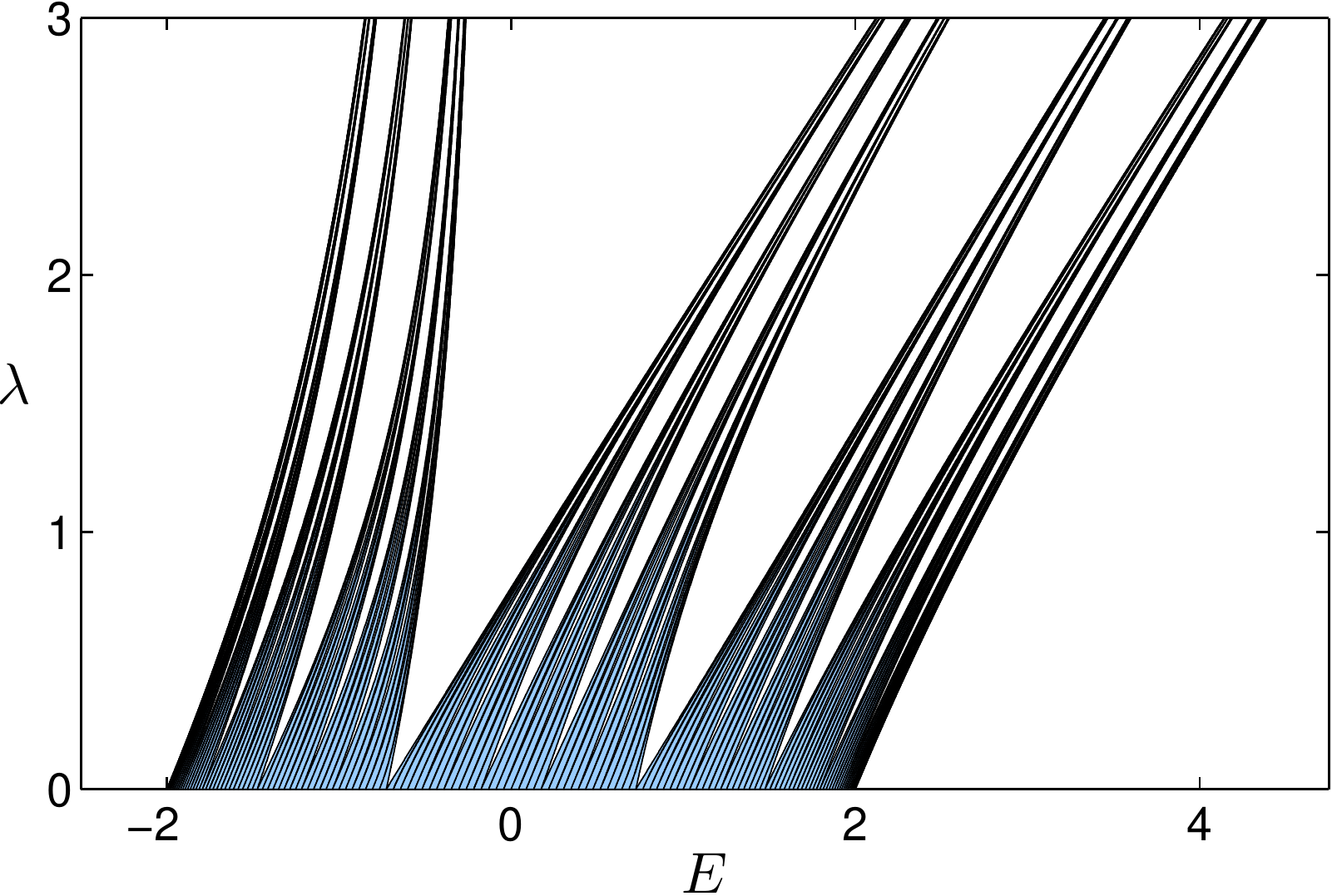}\bigskip}
\centerline{\includegraphics[width=0.8\columnwidth]{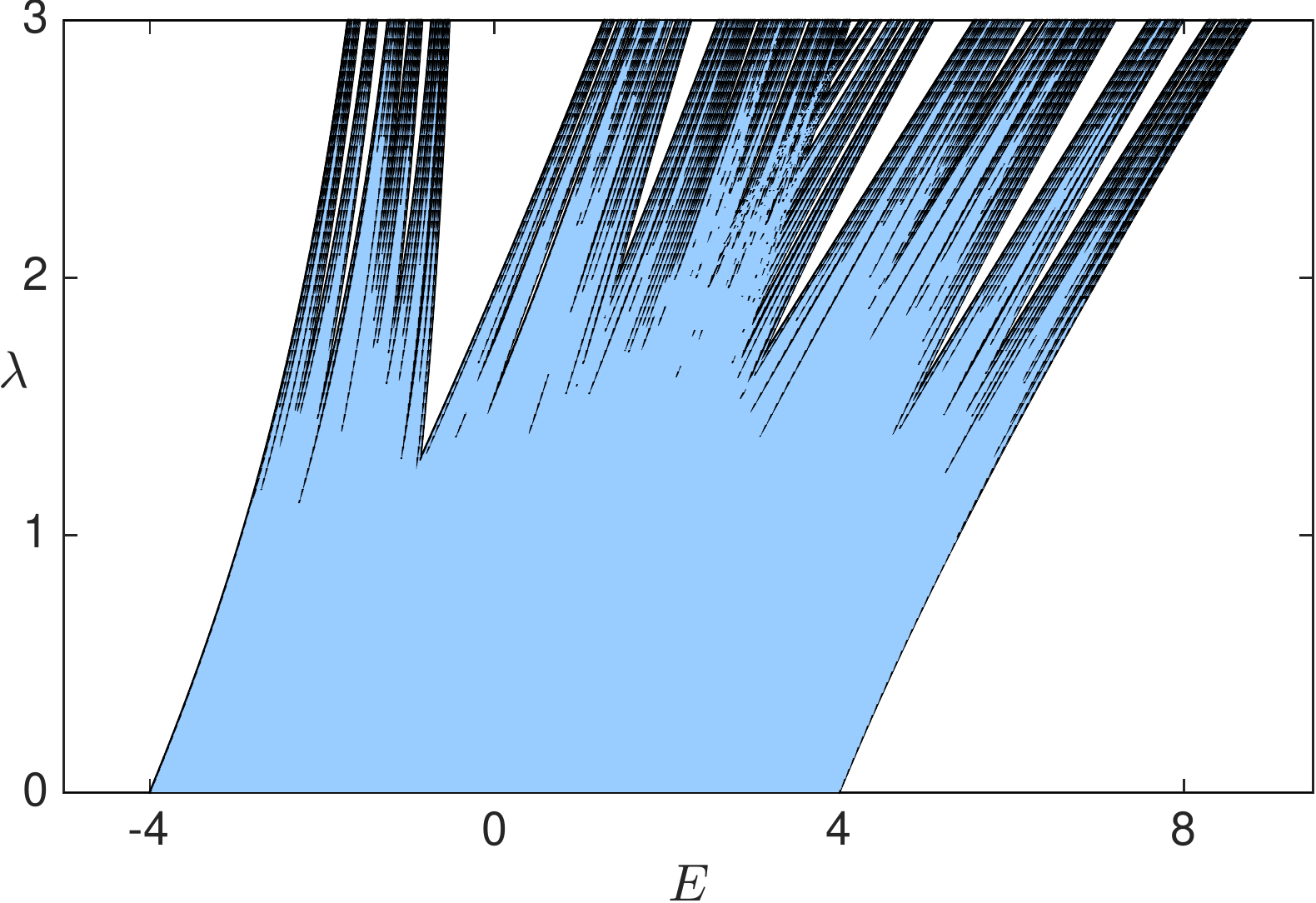}}
\caption{Numerical approximation of the spectrum of the
Fibonacci Hamiltonian (top) and the 2D model obtained as the Cartesian
product of two 1D Fibonacci Hamiltonians (bottom). The $x$-axis corresponds
to the energy $E$, while the $y$-axis corresponds to the coupling constant
$\lambda$. The plots illustrate the instant opening of a dense set of
gaps for the 1D model as the potential is turned on, whereas for the 2D
model there are no gaps in the spectrum for all sufficiently small
$\lambda$. Images courtesy of Mark Embree.\label{fig:fibo}}
\end{figure}

It is a somewhat surprising insight that the \emph{diffraction
  measure}, which is designed to reveal as much as possible about the
distributions of points and is thus clearly not invariant under
topological conjugacy, and the \emph{dynamical spectrum}, which is an
important invariant under metric conjugation of dynamical systems and
thus blind to details of the representative chosen, have such an
important `overlap'.  Consequently, one can translate various results
from either point of view to the other, and profit from this
connection.

While the relation between dynamical and diffraction spectra is by now
well understood, it continues to be an intriguing open problem to find
the connection between these spectra and the spectra of
Schr\"{o}dinger operators associated with aperiodic structures.  The
interest in the latter arises from quantum transport questions in
aperiodically ordered solids; see the article by Damanik, Embree and
Gorodetski in \cite{KLS} for a survey.  In particular, anomalous
transport properties have for a long time been expected (and are
observed in experiments), and could recently be rigorously confirmed
in simple one-dimensional models. Concretely, on the Hilbert space
$\ell^{}_{2} (\mathbb{Z})$, let us consider the Fibonacci Hamiltonian
\[
    (H\psi)_{n} \, = \, \psi_{n-1} + \psi_{n+1} + \lambda v_{n} \psi_{n}
\]
with potential $v_{n}=\chi^{}_{[1-\alpha,1)}(n\alpha\bmod 1)$ (with
constant $\alpha=(\sqrt{5}-1)/2$, the inverse of the golden ratio),
which alternatively could be generated by the Fibonacci substitution
$0\mapsto 1$, $1\mapsto 10$.  For $\lambda>8$, quantum states display
anomalous transport in the sense that they do not move ballistically
or diffusively, nor do they remain localized \cite{DT}. All spectral
measures associated with $H$ are purely singular continuous, while
diffraction and dynamical spectrum are pure point in this case.

Quantitative results about the local and global Hausdorff dimension of
the spectrum and the density of states measure are now available for
all values of the coupling constant $\lambda$. On the other hand,
similar results are currently entirely out of reach for
Schr\"{o}dinger operators associated with the Penrose tiling. However,
there is recent progress for higher-dimensional models obtained as a
Cartesian product of Fibonacci Hamiltonians; see Figure~\ref{fig:fibo}
for an illustration, in particular (of proven properties of the
spectrum) for two dimensions, and \cite{DGS} as well as the article by
Damanik, Embree and Gorodetski in \cite{KLS} for details and further
references.

What we have sketched here is just one snapshop of a field with many
facets and new developments, as evidenced by the contributions to
\cite{Nato,BM,PF,KLS}. Connections exist with many branches of
mathematics, including discrete geometry \cite{CBG}, topology
\cite{Sad} and ergodic theory \cite{Rob,Sol}, to name but a few.
Aperiodic order thus provides a versatile platform for cooperations
and proves the point that mathematics is a unit, and not a collection
of disjoint disciplines.

\begin{small}

\end{small}

\end{document}